\documentclass[11pt,english,a4paper]{smfart}

\usepackage[french,main=english]{babel}
\usepackage[utf8]{inputenc}

\usepackage{graphicx}
\usepackage{ae,amsfonts,euscript,enumerate}
\usepackage{amssymb}
\usepackage{amsmath}
\usepackage[cm]{aeguill}
\usepackage{xcolor}

\renewcommand{\theequation}{\thesection.\arabic{equation}}

\newcommand{\bZ}{\mathbb{Z}}
\newcommand{\bQ}{\mathbb{Q}}

\newcommand{\bC}{\mathbb{C}}
\newcommand{\cK}{\mathcal{K}}

\newtheorem{theo}{Theorem}[section]

\newtheorem{lem}[theo]{Lemma}

\newtheorem{thmx}{Theorem}
\newtheorem{conjx}{Conjecture}

\numberwithin{equation}{section}
\usepackage{fullpage}
\usepackage[T1]{fontenc}

\usepackage{a4wide}

\author{\'E. Delaygue}
\address{
Universit\'e Claude Bernard Lyon 1, CNRS, Centrale Lyon, INSA Lyon, Universit\'e Jean Monnet, ICJ UMR5208, 69622 Villeurbanne, France}
\email{delaygue@math.univ-lyon1.fr}

\title{On Ruzsa's conjecture on congruence preserving functions}
\date{}



\begin{document}

\begin{abstract}
Ruzsa's conjecture asserts that any sequence $(a_n)_{n \geq 0}$ of integers that preserves congruences, \textit{i.e.}, satisfies $ a_{n+k} \equiv a_n \mod k $, and has the growth condition $\limsup_{n \to +\infty} |a_n|^{1/n} < e$, must be a polynomial sequence. While previous results by Hall, Ruzsa, Perelli, and Zannier have confirmed this conjecture under stricter growth bounds, the general case remains open. In this paper, we establish a new partial result by proving that if in addition the generating series $ f = \sum_{n \geq 0} a_n x^n $ has at most two singular directions at $ x = 0 $, then $(a_n)_{n \geq 0}$ is necessarily a polynomial sequence. Our approach is based on an adaptation of Carlson's method, originally developed for the P\'olya--Carlson dichotomy, combined with a refined analysis of Hankel determinants. Specifically, we derive an upper bound on these determinants using P\'olya's inequality and a transfinite diameter argument of Dubinin, while a non-Archimedean divisibility condition on Hankel determinants yields a lower bound, ultimately leading to the rationality of $ f $. This confirms that counterexamples to Ruzsa's conjecture, if they exist, must exhibit at least three singular directions.
\end{abstract}

\maketitle

%%%%%%%%%%
\section{Introduction}
%%%%%%%%%%

Let $(a_n)_{n\geq 0}$ be a sequence of integers that preserves congruences, that is, for all natural integers $n$ and $k$, we have
$$
a_{n+k}\equiv a_n\mod k.
$$
Such a sequence is also known as a \textit{pseudo-polynomial} \cite{Hal71} since, if $P\in\bZ[x]$, then the sequence of integers $(P(n))_{n\geq 0}$ preserves congruences. We shall say that $(P(n))_{n\geq 0}$ is a \textit{polynomial} sequence. A famous conjecture of Ruzsa says that there is a growth condition that forces pseudo-polynomials to be polynomial sequences.

\begin{conjx}[Ruzsa \cite{Ruz71}]
Let $(a_n)_{n\geq 0}$ be a pseudo-polynomial such that
$$
\underset{n\to +\infty}{\limsup}\,|a_n|^{1/n} < e.
$$
Then $(a_n)_{n\geq 0}$ is a polynomial sequence.
\end{conjx}

This conjecture is sharp in the sense that Hall \cite[p. 76]{Hal71} sketched an inductive construction of a pseudo-polynomial $(a_n)_{n\geq 0}$ which is not a polynomial sequence but satisfies
$$
\underset{n\to +\infty}{\limsup}\, |a_n|^{1/n}\leq e.
$$

Ruzsa's conjecture has been proven under more restrictive growth assumptions. In 1971, Hall \cite{Hal71} and Ruzsa \cite{Ruz71} independently established this result when the bound $e$ is replaced by $e-1$. In 1984, Perelli and Zannier demonstrated that the conjecture holds if $e$ is replaced by $e^{0.66}$. Finally, in 1996, Zannier further improved this bound to $e^{0.75}$. To the best of our knowledge, no new cases of this conjecture have been proven since.
\medskip

The results of Perelli and Zannier actually hold in a slightly more general setting, where congruences are considered only modulo prime numbers. Following \cite{DR22}, we define a sequence of integers $(a_n)_{n\geq 0}$ as a \textit{primary pseudo-polynomial} if, for every natural number $n$ and every prime $p$, we have
$$
a_{n+p}\equiv a_n\mod p.
$$
Perelli and Zannier's approach first focused on proving the following elegant result.

\begin{thmx}[Perelli--Zannier \cite{PZ84}]\label{theo: PZ}
Let $(a_n)_{n\geq 0}$ be a primary pseudo-polynomial such that
$$
\underset{n\to +\infty}{\limsup}\,|a_n|^{1/n} < e.
$$
Then $(a_n)_{n\geq 0}$ is $P$-recursive.
\end{thmx}

We recall that a sequence $(a_n)_{n\geq 0}$ is said to be \textit{$P$-recursive} if there exist $d\in\bZ_{\geq 0}$ and polynomials $P_0, \dots, P_d$ in $\mathbb{Q}[x]$, not all zero, such that, for every natural number $n$, we have
$$
P_d(n)\,a_{n+d}+\cdots + P_0(n)\,a_n=0. 
$$
Equivalently \cite{Sta80}, this means that the generating series $f$ of $(a_n)_{n\geq 0}$ is annihilated by a nonzero linear differential operator with coefficients in $\mathbb{Q}[x]$, or, equivalently, that $f$ is \textit{$D$-finite}.

Thus, under the assumptions of Theorem \ref{theo: PZ}, the generating series $f$ of $(a_n)_{n\ge 0}$
is $D$-finite and has a strictly positive radius of convergence at the origin. It follows that $f$ satisfies a
non-trivial linear differential equation
$$
p_r(x) f^{(r)}(x) + p_{r-1}(x) f^{(r-1)}(x) + \cdots + p_0(x) f(x)=0,
$$
with $p_i(x)\in\bQ[x]$, $0\leq i\leq r$, and $p_r\not\equiv 0$. Let $S:=\{x\in\bC:\ p_r(x)=0\}$. A standard continuation theorem for linear differential equations with analytic coefficients implies that $f$ admits analytic continuation along any path in $\bC\setminus S$,  see \emph{e.g.} \cite[Ch.~5]{vdPS03}. 

We define $\mathrm{Sing}(f)$ to be the set of finite singularities of $f$, \emph{i.e.} the set of points $\omega\in\bC$ such that no analytic continuation of $f$ from $0$ extends holomorphically to a neighbourhood of $\omega$. Then $\mathrm{Sing}(f)$ is a subset of $S$ and has at most $\deg p_r$ elements. We call a \emph{singular direction} for $f$
at $0$ the principal argument of any $\omega\in\mathrm{Sing}(f)$. The aim of this paper is
to establish the following new partial case of Ruzsa's conjecture.

\begin{theo}\label{theo: main}
Let $(a_n)_{n\geq 0}$ be a primary pseudo-polynomial such that
$$
\underset{n\to +\infty}{\limsup}\,|a_n|^{1/n} < e.
$$
If its generating series has at most two singular directions at $0$, then $(a_n)_{n\geq 0}$ is a polynomial sequence.
\end{theo}

Thus, a counterexample to Ruzsa's conjecture could only exist if it has at least three singular directions at $0$.
Perelli and Zannier proved in Theorem \ref{theo: PZ} that the generating series $f$ of $(a_n)_{n\geq 0}$ is $D$-finite, but their argument is non-effective and does not provide any explicit control on the size of $\mathrm{Sing}(f)$.
An effective refinement obtained in joint work with Rivoal~\cite[Theorem~3]{DR22} yields explicit bounds for the order and the degrees of an annihilating differential operator, hence an explicit upper bound for the size of $\mathrm{Sing}(f)$.
But this bound is typically very large when $\limsup_{n\to\infty}|a_n|^{1/n}$ is close to $e$. In particular, this does not allow one to apply Theorem~\ref{theo: main} without the restriction to at most two singular directions.

Our approach to proving Theorem \ref{theo: main} is based on Carlson's method \cite{Car21} for the famous P\'olya--Carlson dichotomy (see Theorem \ref{theo: PC} below), a method that was recently used by Dimitrov \cite{Dim19} to prove the Schinzel--Zassenhaus conjecture on polynomials. 

Zannier's approach \cite{Zan96} to establishing Ruzsa's conjecture with the bound $e^{0.75}$ instead of $e$ is of a different nature; it relies on deep results on the arithmetic of linear differential equations.

%%%%%%%%%%
\section{A P\'olya--Carlson approach}
%%%%%%%%%%

Under the hypotheses of Theorem \ref{theo: main}, we shall show that the generating series $f$ of $(a_n)_{n\geq 0}$ is the power series expansion of a rational fraction. Using the congruences satisfied by $(a_n)_{n\geq 0}$, it then becomes straightforward to prove that $f$ is the power series expansion of a rational fraction whose denominator is a power of $(1-x)$. This implies that $(a_n)_{n\geq 0}$ is a polynomial sequence, as desired. In fact, Zannier showed in \cite[p. 396]{Zan96} that if the generating series $f$ of a primary pseudo-polynomial $(a_n)_{n\geq 0}$ is algebraic over $\mathbb{Q}(x)$, then $(a_n)_{n\geq 0}$ is a polynomial sequence.

%%%%%%%%%%
\subsection{Kronecker's rationality criterion}
%%%%%%%%%%

Let $f=\sum_{n=0}^\infty a_nx^n\in\bC[[x]]$ be a formal power series with complex coefficients. For every positive integer $n$ we denote by $H_n(f)$ the $n$-th \textit{Hankel matrix} of $f$ defined by
$$
H_n(f)=\left(\begin{array}{cccc}
a_0 & a_1 & \cdots & a_{n-1}\\
a_1 & a_2 & \cdots & a_n	\\
\vdots & \vdots & & \vdots \\
a_{n-1} & a_n & \cdots & a_{2n-2}
\end{array}\right).
$$
We shall say that $\det H_n(f)$ is the $n$-th \textit{Hankel determinant} of $f$. Kronecker proved in \cite{Kro81} the following.

\begin{thmx}[Kronecker's rationality criterion]\label{theo: Kronecker}
A power series $f\in\bC[[x]]$ is the expansion of a rational fraction if and only if almost all Hankel determinants of $f$ are zero.
\end{thmx}

When $f\in\bZ[[x]]$, all Hankel determinants of $f$ belong to $\bZ$, which provides non-archimedean bounds: for every nonnegative integer $n$ and every prime $p$, we have 
$$
|\det H_n(f)|_p\leq 1, 
$$
where $|\cdot|_p$ stands for the usual $p$-adic absolute value on $\bQ$. If one manages to prove that $|\det H_n(f)| < 1$ for large enough $n$, then we have $\det H_n(f)=0$ for such $n$ and $f$ is the power series expansion of a rational fraction. This archimedean upper bound can be achieved \textit{via} P\'olya's inequality under additional assumptions. This was Carlson's idea \cite{Car21} to prove the following dichotomy.

\begin{thmx}[P\'olya--Carlson dichotomy]\label{theo: PC}
A power series $f\in\bZ[[x]]$ that converges inside the unit disk is either rational or admits the unit circle as a natural boundary. 
\end{thmx}

Our approach to proving Theorem \ref{theo: main} is similar. We will first establish an archimedean upper bound on the Hankel determinants of $f$ using P\'olya's inequality. We will then derive a non-archimedean bound using the congruences satisfied by $(a_n)_{n\geq 0}$, which will allow us to conclude.

%%%%%%%%%%
\subsection{An archimedean upper bound}
%%%%%%%%%%

The aim of this section is to prove the following archimedean upper bound.

\begin{lem}\label{lem: ArchiBound}
Let $(a_n)_{n\geq 0}$ be a sequence of complex numbers such that its generating series $f$ is $D$-finite with a positive radius of convergence $\rho$ and at most $r$ singular directions at $0$. Then, we have
$$
\underset{n\rightarrow +\infty}{\limsup}\,|\det H_n(f)|^{1/n^2}\leq \frac{1}{4^{1/r}\rho}.
$$
\end{lem}

The proof of the P\'olya--Carlson dichotomy and of Lemma \ref{lem: ArchiBound} is based on the following well-known inequality.

\begin{thmx}[P\'olya's inequality \cite{Pol28}]\label{theo: Polya}
Let $(a_n)_{n\geq 0}$ be a sequence of complex numbers and $f$ its generating series. Let $K$ be a compact subset of $\bC$ such that $\bC\setminus K$ is connected and 
$$
g(x)=\sum_{n=0}^\infty a_n x^{-n-1}
$$
is analytic on $\bC\setminus K$. Then, we have
$$
\underset{n\rightarrow +\infty}{\limsup}\,|\det H_n(f)|^{1/n^2}\leq \tau(K),
$$
where $\tau(K)$ is the transfinite diameter of $K$.
\end{thmx}

In our context, we follow the approach of Dimitrov \cite{Dim19} and will use the following theorem of Dubinin to estimate the transfinite diameter of a hedgehog
$$
\cK(a_1,\dots,a_r)=\bigcup_{i=1}^r\left[0,a_i\right],
$$
where $a_1,\dots,a_r$ are complex numbers.

\begin{thmx}[Dubinin \cite{Dub85}]\label{theo: Dubinin}
The hedgehog $\cK(a_1,\dots,a_r)\subset\bC$ has transfinite diameter at most $(\max_i |a_i|^r/4)^{1/r}$, with equality if and only if the points $a_1,\dots,a_r$ form the vertices of a regular $r$-gon centered at the origin.
\end{thmx}

\begin{proof}[Proof of Lemma \ref{lem: ArchiBound}]
Since $f$ is $D$-finite with a positive radius of convergence $\rho$ and at most $r$ singular directions, it can be analytically continued to the entire complex plane, except along a finite number of pairwise disjoint radial half-lines originating from some singularities $\alpha_1, \dots, \alpha_r$ satisfying
$$
\min_i\,|\alpha_i|=\rho.
$$
Making the change of variables $x\mapsto 1/x$, we obtain that
$$
\frac{1}{x}f\left(\frac{1}{x}\right)=\sum_{n=0}^\infty a_n x^{-n-1}
$$
can be analytically continued to $\bC\setminus\cK$ where $\cK$ is the hedgehog 
$$
\cK\left(\frac{1}{\alpha_1},\dots,\frac{1}{\alpha_r}\right)=\bigcup_{i=1}^r\left[0,\frac{1}{\alpha_i}\right].
$$
We can then use Theorems \ref{theo: Polya} and \ref{theo: Dubinin} to obtain the upper bound 
$$
    \underset{n\rightarrow +\infty}{\limsup}\,|\det H_n(f)|^{1/n^2}\leq \frac{1}{4^{1/r}\rho},
$$
as desired.
\end{proof}

%%%%%%%%%%
\subsection{A non-archimedean bound}
%%%%%%%%%%

The aim of this section is to prove the following non-archimedean bound.

\begin{lem}\label{lem: non-ArchiBound}
Let $(a_n)_{n\geq 0}$ be a primary pseudo-polynomial and $f$ its generating series. Then, for every positive integer $n$, $\det H_n(f)$ is divisible by
$$
\prod_{\substack{p \le n-1 \\ p \text{ prime}}} p^{n - p}.
$$
\end{lem}

The main idea behind Lemma \ref{lem: non-ArchiBound} is to consider the binomial transform \( g \) of \( f \), for which it is easier to obtain large divisibility properties of its Hankel determinants. If $(a_n)_{n\geq 0}$ is a sequence of complex numbers, then we consider its \textit{binomial transform} $(b_n)_{n\geq 0}$ defined by
$$
b_n:=\sum_{k=0}^n(-1)^{n-k}\binom{n}{k}a_k.
$$
One can also recover $(a_n)_{n\geq 0}$ from $(b_n)_{n\geq 0}$ \textit{via} the formula
$$
a_n=\sum_{k=0}^n\binom{n}{k}b_k.
$$

In the remainder of the article, if $f$ denotes the generating series of a sequence $(a_n)_{n\geq 0}$ of complex numbers and $g$ represents the generating series of the binomial transform $(b_n)_{n\geq 0}$ of $(a_n)_{n\geq 0}$, then we shall refer to $g$ as the \textit{binomial transform} of $f$. It is well-known that, for every positive integer $n$, we have $\det H_n(f)=\det H_n(g)$, see \emph{e.g.} \cite[Lemma~15]{Kra99}. 

\begin{proof}[Proof of Lemma \ref{lem: non-ArchiBound}]
Let $(a_n)_{n\geq 0}$ be a primary pseudo-polynomial, $(b_n)_{n\geq 0}$ its binomial transform, and let $f$ and $g$ be their respective generating series. According to \cite[Theorem 1]{DR22}, for every natural integer $n$, $b_n$ is divisible by the product $P_n$ of the primes less than or equal to $n$.  

Since for every positive integer $n$, $P_n$ divides $b_n$, the $i$-th row of $H_n(g)$ can be factored by $P_{i-1}$. It follows that $\det H_n(g)$ is divisible by 
$$
\prod_{k=1}^{n-1}P_k = \prod_{\substack{p \le n-1 \\ p \text{ prime}}} p^{n-p},
$$
and so is $\det H_n(f)$. 
\end{proof}

We can now prove our main result.

%%%%%%%%%%
\subsection{Proof of Theorem \ref{theo: main}}
%%%%%%%%%%

Let $(a_n)_{n\geq 0}$ be a primary pseudo-polynomial such that 
$$
\underset{n\to+\infty}{\limsup}\,|a_n|^{1/n} < e,
$$
and let $f$ be its generating series. According to the work of Perelli and Zannier \cite{PZ84}, $f$ is $D$-finite with a radius of convergence $\rho> 1/e$. Since $f$ has at most two singular directions at $0$, we can apply Lemma \ref{lem: ArchiBound} and obtain that
\begin{equation}\label{eq: RealArchiBound}
\underset{n\rightarrow +\infty}{\limsup}\,|\det H_n(f)|^{1/n^2}< \frac{e}{2}.
\end{equation}
Furthermore, according to Lemma \ref{lem: non-ArchiBound}, for all $n$ in $\mathbb{N}$, $\det H_n(f)$ is divisible by
\begin{equation}\label{eq: prodP}
\prod_{\substack{p \le n-1 \\ p \text{ prime}}} p^{n-p}.
\end{equation}
It therefore remains to determine an asymptotic estimate of the quantity \eqref{eq: prodP}. Taking the natural logarithm, we get
\begin{align*}
\log\,\prod_{p\leq n-1}p^{n-p}
&=\sum_{p \leq n-1} (n - p)\log p\\
&=n\sum_{p \le n-1} \log p 
- \sum_{p \le n-1} p\log p.
\end{align*}
Let 
$$
\theta(x) = \sum_{p \le x} \log p
$$
be the first Chebyshev function. Then
$$
n \sum_{p \le n-1} \log p 
= n\,\theta(n-1).
$$
To handle the term $\sum_{p \le n-1} p\log p$, we use Abel transformation. Let $(c_n)_{n\geq 0}$ be given by $c_n=\log n$ if $n$ is prime and $c_n=0$ otherwise. Then we have
\begin{align*}
\sum_{p \le n-1} p\log p
&=\sum_{k=0}^{n-1}k\,c_k\\
&= n\sum_{k=0}^{n-1}c_k-\sum_{k=0}^{n-1}\left(\sum_{i=0}^k c_i\right)(k+1-k)\\
&= n\theta(n-1) - \sum_{k=0}^{n-1}\theta(k),
\end{align*}
leading to
$$
\log\,\prod_{p\leq n-1}p^{n-p}
= \sum_{k=0}^{n-1}\theta(k).
$$
The Prime Number Theorem states that $\theta(x) \sim x $ as $x \to +\infty$. Hence, as $n\to +\infty$, we have
$$
\sum_{k=0}^{n-1}\theta(k)\sim\frac{n^2}{2},
$$
and
\begin{equation}\label{eq: asympto}
\prod_{p\leq n-1}p^{n-p} = \exp\left(\frac{n^2}{2}+o(n^2)\right).
\end{equation}

Since $\sqrt{e} > e/2$, the bounds \eqref{eq: RealArchiBound} and \eqref{eq: asympto} force $\det H_n(f)$ to be zero for sufficiently large $n$. By Kronecker's rationality criterion, Theorem \ref{theo: Kronecker}, $f$ is the series expansion of a rational fraction. As mentioned above, Zannier showed in \cite[p. 396]{Zan96} that if the generating series $f$ of a primary pseudo-polynomial $(a_n)_{n\geq 0}$ is algebraic over $\mathbb{Q}(x)$, then $(a_n)_{n\geq 0}$ is a polynomial sequence, which completes the proof of Theorem \ref{theo: main}.
$\hfill\square$
\medskip

\noindent\textbf{Acknowledgements}. We warmly thank Dang-Khoa Nguyen for bringing Dubinin's result to our attention and Jiang Zeng for pointing us to \cite[Lemma~15]{Kra99}.

\end{document}